\documentclass[12pt]{article}

\usepackage{amsmath}
\usepackage{amsfonts}
\usepackage{amsthm}

\numberwithin{equation}{section}

\def\Cmpx{{\mathbb{C}}}
\def\Real{{\mathbb{R}}}
\def\Intg{{\mathbb{Z}}}
\def\cnj{\overline}
\def\union{{\cup}}

\def\tfrac#1#2{{\textstyle{\frac{#1}{#2}}}}

\def\scrhalf{{\raisebox{.3ex}%
{$\scriptstyle 1$}\!/\!\raisebox{-.3ex}{$\scriptstyle 2$}}}

\def\wt{\widetilde}
\def\wh{\widehat}
\def\Im{{\mathrm{Im}}}

\def\spanrm{{\mathrm{span}}}
\def\innerprod(#1){{\langle #1 \rangle}}
\def\norm#1{\|#1\|}
\def\rvec(#1){{|#1\rangle}} 
\def\lvec(#1){{\langle#1|}}

\newtheorem{theorem}{Theorem}

\def\ad#1{{{\mathrm{ad}}_{#1}}}
\def\Pexpr{{\cal P}}
\def\Qexpr{{\cal Q}}
\def\Rp{{\wh R}}
\def\Alg{{\cal A}}
\def\Man{{\cal M}}
\def\Sph{{S^2}}
\def\Hypone{{{\cal H}_1}}
\def\Hyptwop{{{\cal H}_2^+}}
\def\Hyptwom{{{\cal H}_2^-}}
\def\PSph{{\Pexpr_\Sph}}
\def\PHone{{\Pexpr_\Hypone}}
\def\PHtwo{{\Pexpr_{{\cal H}_2^\pm}}}
\def\Wigner(#1,#2,#3){{\left<
\begin{array}{ccc} & #2 & \\ #1 & & 0 \\ & #3 & \end{array}\right>}}
\def\CG(#1,#2,#3,#4,#5,#6){{C^{#1}_{#4}{}^{#2}_{#5}{}^{#3}_{#6}}}
\def\RMP(#1,#2,#3){{{\cal R}^{#1}{}^{#2}{}^{#3}}}
\def\WsixJ(#1,#2,#3,#4,#5,#6)%
{\left\{ \begin{array}{ccc} #1\ &\ #2\ &\ #3
\\#4\ &\ #5\ &\ #6\end{array}\right\}}

\title{``Wick Rotations'': The Noncommutative Hyperboloids and Other
Surfaces of Rotations}

\author{{\bf Jonathan  Gratus%
\thanks{Funded by the Royal Society of London European 
Science Exchange Programme}}%
\\
Laboratoire de Gravitation et Cosmologie Relativistes%
\thanks{\it Laboratoire associ\'e au CNRS {\rm URA 769}}
\\
Tour 22/12 4eme etage, Boite Courrier142, 4pl Jussieu. F75252
Paris%
\thanks{Present address: Department of Physics, Lancaster University,
Lancaster LA1 4YB, UK.}
\\
email: jg@luna.ph.lancs.ac.uk
}
\date{January 8, 1998}

\begin{document}

\maketitle

\begin{abstract}
A ``Wick rotation'' is applied to the noncommutative sphere to produce
a noncommutative version of the hyperboloids. A harmonic basis of the
associated algebra is given.  It is noted that, for the one sheeted
hyperboloid, the vector space for the noncommutative algebra can be
completed to a Hilbert space, where multiplication is not continuous.
A method of constructing noncommutative analogues of surfaces of
rotation, examples of which include the paraboloid and the
$q$-deformed sphere, is given. Also given are mappings between
noncommutative surfaces, stereographic projections to the complex
plane and unitary representations.  A relationship with one
dimensional crystals is highlighted.

\end{abstract}

\tableofcontents


\section{Introduction}

This letter is divided into two sections. The first is concerned with
analytically continuing the algebra of the noncommutative sphere so
producing the noncommutative analogue of the hyperboloid, whilst the
second section uses this to produce noncommutative analogues of a vast
collection of axially symmetric two dimensional surfaces.

\vskip 1em

As every school child knows $x^2+y^2+z^2=R^2$ is the equation for a
sphere ($S^2$) of radius $R$ embedded in $\Real^3$. Likewise
$z^2-x^2-y^2=R^2$ is the equation of a two sheeted hyperboloid
($\Hyptwop\union\Hyptwom$ where $\Hyptwop$ and $\Hyptwom$ are the
upper and lower sheets), and $z^2-x^2-y^2=-R^2$ is the equation for
the one sheeted hyperboloid $\Hypone$. It is obvious that if one
performs the ``Wick rotation'' $x\to ix$ and $y\to iy$ one passes from
the sphere to the two sheeted hyperboloid, whilst the substitution
$R\to iR$ takes one from the two sheeted hyperboloid to the one
sheeted hyperboloid. 

The standard method of analysing the noncommutative or ``fuzzy''
sphere is by the use of matrices~\cite[chapter 7.2]{Madore_book}. In
such an approach it is not clear how one can perform a ``Wick
rotation''.  However in \cite{Gratus5}, we present a two parameter
algebra $\Pexpr(\varepsilon,R)$ which may be thought of as the
noncommutative sphere, since for $\varepsilon=0$, $\Pexpr(0,R)$ is
equivalent to the algebra of complex valued functions on the sphere.
For a discreet set of $\varepsilon$, $\Pexpr(\varepsilon,R)$ can be
mapped into the algebra of matrices. Since the approach of that
article is more algebraic it is easier to perform a ``Wick Rotation''
so producing noncommutative analogues of the one and two sheeted
hyperboloid.  In section \ref{ch_Wick} we give the details of such a
rotation. The new algebra contains an extra parameter,
$\alpha\in\Cmpx$ with $|\alpha|=1$ which gives the angle of rotation,
smoothly rotating between the algebra for the sphere ($\PSph$ when
$\alpha=1$) and the algebra for the hyperboloids ($\PHtwo$ and
$\PHone$ when $\alpha=i$). We rewrite the major expressions in
\cite{Gratus5} for general $\alpha$.  We also give a formula for the
product of two basis polynomials in terms of Wigner $6j$ symbols.

The one sheeted hyperboloid $\Hypone$ is of particular interest to
physicists since it may be considered as a globally hyperbolic
spacetime in one plus one dimensions and as a two dimensional
equivalent of de Sitter space. The algebra $\PHone$ associated with
this space may aid the construction of a noncommutative (quantum)
theory of fields on de Sitter spaces. This algebra has a very useful
property. The sesquilinear form on $\PHone$ is positive definite and
hence an inner product. It is therefore possible to complete the
underlying vector space to produce a Hilbert space $\cnj\PHone$.
Multiplication within $\cnj\PHone$ is not continuous and, as a result,
the elements of $\PHone$ may be represented by unbounded operators as
they act on $\cnj\PHone$ by left (or right) multiplication.  We
discuss the existence or otherwise of a representation of $su(1,1)$ by
the action of left multiplication on $\cnj\PHone$.

\vskip 1em

In the second section of this letter we construct noncommutating
analogies for surfaces of rotation. Connected surfaces of rotation are
either topologically equivalent to the sphere, the disc or the
cylinder.

In section \ref{ch_suf} we give a definition of noncommutative
surfaces of rotation and show how to map functions on one surface to
functions on another.  These maps are generalisations of the
Holstein and Primakoff formalism.  They indicate a strong
relationship between (1) the topology of the manifold (2) the
Hermitian conjugation of the algebra, and (3) the unitary
representations of the algebra.

In section \ref{ch_suf_eg} we show that the algebras for the
noncommutative sphere and hyperboloids analysed in section
\ref{ch_Wick} fit nicely into this framework and that the
Heisenberg-Weil algebra can be viewed as the noncommutative
paraboloid. We also show that the $q$-deformed sphere may be seen as a
way of continuously deforming the sphere (for $q=1$) into a cylinder
with end discs (for $q=\infty$).

In section \ref{ch_suf_rep} we give the unitary representations of the
noncommutative surfaces. Compact surfaces have finite dimensional
representations, whilst non compact surfaces have infinite dimensional
representations. We highlight a relationship between representations
of noncommutative surfaces and a large class of one dimensional
crystal lattice problems.

In section \ref{ch_suf_Cmplx} we show that the algebra for a
noncommutative surface may also be regarded as the algebra for a
noncommutative complex variable. That is, in certain situations, we
can construct a one parameter algebra in which a complex variable $z$
does not commute with is complex conjugate $\cnj z$. In the
commutative limit this algebra is equivalent to the algebra of complex
valued functions on a domain of $\Cmpx$. We give explicit maps between
$\PSph$ and $\PHtwo$ and the noncommutative complex plane. These maps
may be regarded as the noncommutative analogue of stereographic
projections.

\vskip 1em

Finally in section \ref{ch_Disc}, we discuss how one might use the
results in this letter to develop a quantum theory of gravity, some of
the problems that are likely to arise, and some of their possible
solutions.


\section{A ``Wick rotation'' of the Noncommutative Sphere to the
Noncommutative Hyperboloid}
\label{ch_Wick}

Let us write $J_0=z$ and $J_\pm=x\pm
i y $ then $J_0^2+\tfrac12 (J_+J_-+J_-J_+)=R^2$ is the
equation of a sphere, whilst $J_0^2-\tfrac12 (J_+J_-+J_-J-)=\pm R^2$
are the equations of the hyperboloids. The ``Wick rotation'' from a
sphere to the hyperboloids may be made by mapping $J_\pm\mapsto
iJ_\pm$ and allowing $R$ to be complex. This mapping may also be
continued to the noncommutative case. To make the rotation more
explicit we consider $J_\pm\mapsto \alpha J_\pm$ where $\alpha$ may be
any complex nonzero number. However, since we can rescale $J_\pm$ we
shall set $|\alpha|=1$.

For the case of the noncommutative sphere there exists an algebra of
polynomials $\Pexpr$ given in \cite{Gratus5}. This algebra now becomes
the algebra of polynomials generated by $\{J_+,J_-,J_0\}$ where
\begin{align}
[J_0,J_+] &= \varepsilon J_+ &
[J_0,J_-] &= -\varepsilon J_- &
[J_+,J_-] &= 2\varepsilon\alpha^2 J_0 &
J_0^2 + \frac1{2\alpha^2}(J_+J_- + J_-J_+) &= R^2
\label{Wick_com_rel}
\end{align}
The algebra $\Pexpr$ thus depends on $\varepsilon,R,\alpha\in\Cmpx$
which are all independent.

The only way this algebra is distinguished from a simple
complexification of the case when $\alpha=1$ is by the choice of
Hermitian conjugate.  This is given by $\dagger:\Pexpr\mapsto\Pexpr$
\begin{align}
J_0^\dagger &=J_0, &
J_+^\dagger &=J_-, & 
J_-^\dagger &=J_+, &
(ab)^\dagger &= b^\dagger a^\dagger &
\lambda^\dagger &= \cnj{\lambda} &
\forall\, a,b\in\Pexpr,\, \lambda\in\Cmpx
\end{align}
Clearly this  conjugation is consistent with (\ref{Wick_com_rel}) if
and only if $\varepsilon,R^2,\alpha^2\in\Real$. There are six cases
when $\varepsilon,R^2,\alpha^2\in\Real$:

\begin{center}
\begin{tabular}{lcll}
$\alpha^2=1$ & $R^2>0$ & $\Pexpr=\PSph$  & Sphere \\
$\alpha^2=1$ & $R^2=0$ & $\Pexpr$  & Point \\
$\alpha^2=1$ & $R^2<0$ & $\Pexpr$  & No Manifold \\
$\alpha^2=-1$ & $R^2>0$ & $\Pexpr=\PHtwo$  & Two-sheeted Hyperboloid \\
$\alpha^2=-1$ & $R^2=0$ & $\Pexpr$  & Two cones \\
$\alpha^2=-1$ & $R^2<0$ & $\Pexpr=\PHone$  & One-sheeted Hyperboloid \\
\end{tabular}
\end{center}

\noindent
The algebra $\Pexpr$ is still valid even when it does not correspond to a
manifold.  In this letter we shall let $\alpha$, $\varepsilon$ and $R$
be formal, self Hermitian ($\alpha^\dagger=\alpha$ etc.)  parameters in
the centre of $\Pexpr$. Thus we can still do manipulations involving
conjugation without requiring them to be real numbers.

The sesquilinear form is defined in the same way as in \cite{Gratus5},
that is $\innerprod(f,g)=\pi_0(f^\dagger g)$ where $\pi_0(f)$ is the
coefficient of unity when $f$ is written as a formally tracefree
symmetric polynomial. For the commutative sphere ($\alpha^2=1$,
$\varepsilon=0$) this is the standard inner product calculated by
integrating over the sphere; $\innerprod(f,g)=\int_\Sph \cnj{f}g\;
d\mu$.  It is also the trace with respect to (\ref{Pmn_rep_kj}); the
finite dimensional representation of $sl(2,\Cmpx)$.  With respect to
this inner product there is an orthogonal (but unnormalised) basis of
$\Pexpr$ given by $\{P^m_n(\varepsilon,R)\,|\ n,m\in\Intg,n\ge0,|m|\le
n\}$ where
\begin{align}
P^m_n(\varepsilon,R) &=
\alpha^{m-n}\varepsilon^{m-n}
\left(\frac{(n+m)!}{(2n)!\,(n-m)!}\right)^{\scrhalf} 
(\ad{J_-})^{n-m}(J_+{}^n)
\label{Pmn_def_Pmn}
\end{align}

When written as a formally tracefree symmetric polynomial in
$(J_0,J_+,J_-)$, $P^m_n$ is homogeneous of order $n$ and is
independent of $R$ and $\varepsilon$ (but not necessarily $\alpha$).
Each $P^m_n$ is an eigenvector of the operators $\ad{J_0}$ and
$\Delta=\ad{J_0}^2 +
\tfrac1{2\alpha^2}(\ad{J_+}\ad{J_-}+\ad{J_-}\ad{J_+})$:
\begin{align}
\ad{J_0} P^m_n &= \varepsilon m P^m_n 
\label{Pmn_ez_Pmn} 
\\
\Delta P^m_n &= \varepsilon^2 n(n+1) P^m_n 
\label{Pmn_Del_Pmn} 
\end{align}
The ladder operators $\ad{J_+},\ad{J_-}$ increase or decrease $m$: 
\begin{align}
\ad{J_\pm} P^m_n &= \alpha\varepsilon (n\mp m)^{\scrhalf} 
(n\pm m+1)^{\scrhalf} P^{m\pm 1}_n 
\label{Pmn_epm_Pmn}
\end{align}
and the normal of $P^m_n$ is given by
\begin{align}
\norm{P^m_n}^2 &=
\alpha^{2n}
\frac{(n!)^2}{(2n+1)!} 
\prod_{r=1}^n(4R^2+\varepsilon^2(1-r^2))
\label{Pmn_Norm}
\end{align}
If we require $\varepsilon,\alpha^2,R^2\in\Real$ then, in general,
$\norm{P^m_n}^2$ may be positive negative or zero. However, for the
1-sheeted hyperboloid ($\alpha^2=-1$, $4R^2<-\varepsilon^2$)
$\norm{P^m_n}^2 >0$ for all $n$.  This enables us to complete $\PHone$
into a Hilbert space denoted $\cnj{\PHone}$. It is clear that the
action of left or right multiplication by $J_0$ or $J_\pm$ on the
$\cnj\PHone$ are given by unbounded operators. This is
examined at the end of this section.

\vskip 1em

The finite dimensional representation of $sl(2,\Cmpx)$ are given for
$2k\in\Intg$, $k\ge 0$ by
\begin{align}
J_0 \rvec(k,j) &= \varepsilon j \rvec(k,j) 
&
J_\pm \rvec(k,j) &= \alpha\varepsilon (k\mp j)^{\scrhalf} 
(k\pm j+1)^{\scrhalf} \rvec(k,j\pm 1) 
\label{Pmn_rep_kj}
\end{align}
This representation is unitary when $\alpha=1$. It is easy to see that
any other representation is unitary only when
$\alpha^2,\varepsilon,R^2\in\Real$. As a result the only other unitary
representations are the classical unitary representations of
$su(1,1)$.

As before the projection $\pi_0:\Pexpr\mapsto\Cmpx$ is given by the
trace: $\pi_0(f)=\tfrac1{2k+1}\sum_{j=-k}^k\lvec(k,j)f\rvec(k,j)
$. This is used to  calculate (\ref{Pmn_Norm}) using
 $R^2=\varepsilon^2 k(k+1)$.  It is not clear how one can
use the unitary representations of $su(1,1)$ to generate a formula for
$\pi_0(f)$.


\begin{theorem}
\label{thm_RM}
As operators on a Hilbert space, $P^m_n$ can be viewed as a Wigner
operator:
\begin{align}
P^m_n\rvec(k,j)
&=
(-1)^{n} 
\norm{P^m_n}
(2n+1)^\scrhalf 
\Wigner(2n,n,n+m)
\rvec(k,j)
\label{RM_Wig_op}
\end{align}
We can use this to write the formula for the product of two basis
elements in terms of Wigner $6j$ symbols:
\begin{align}
P^{m_1}_{n_1}
P^{m_2}_{n_2}
&=
\sum_{n=|n_1-n_2|}^{n=n_1+n_2}
\CG(n_1,n_2,n,m_1,m_2,m_1+m_2)
\RMP(n_1,n_2,n)
P_{n}^{m_1+m_2} 
\label{RM_def}
\end{align}
where 
$\CG(n_1,n_2,n,m_1,m_2,m_1+m_2)$ is the Clebsh-Gordon coefficient, and
the reduced matrix element $\RMP(n_1,n_2,n)$ is
given by
\begin{align}
\RMP(n_1,n_2,n)
&=
(-1)^{2k+n_1+n_2}
\frac{\norm{P^{m_1}_{n_1}}\norm{P^{m_2}_{n_2}}}
{\norm{P^{m_1+m_2}_{n}}}
(2k+1)^\scrhalf 
(2n_1+1)^\scrhalf (2n_2+1)^\scrhalf 
\WsixJ(k,n_1,k,n_2,k,n) 
\label{RM_res}
\end{align}
where the symbol in the curly brackets is Wigner's 6-$j$ coefficient.
\end{theorem}

\begin{proof}
By application of the Wigner-Eckart theorem we have
\begin{align*}
P^m_n \rvec(k,j) &= 
D_{nk} \CG(k,n,k,j,m,j+m) \rvec(k,j+m)
\end{align*}
where $D_{nk}\in\Cmpx$ is the associated reduced matrix element. To calculate
this put $m=n$. In this case $\CG(k,n,k,j,n,j+n)$ has only one
term. Substituting this into the definition of the Wigner operator
\cite[eqn (3.341)]{Biedenharn1} gives
(\ref{RM_Wig_op}). One then uses the product law given by
\cite[eqn (3.350)]{Biedenharn1}.

Here we define
$\norm{P^m_n}\equiv\alpha^n\big(\alpha^{-2n}\norm{P^m_n}^2\big)^\scrhalf$
which is well defined since $\alpha^{-2n}\norm{P^m_n}^2>0$.
\end{proof}

Because of the defining equations for the algebra
(\ref{Wick_com_rel}), one can use $R^2=\varepsilon^2 k(k+1)$ to remove
$k$ from (\ref{RM_res}) to give an expression for the reduced matrix
element $\RMP(n_1,n_2,n)$ which is a polynomial in $R^2$.

These formulae extend naturally to the algebra of deformed rotation
matrices given in \cite{Gratus6}.


\subsection*{A possible unitary representation of $su(1,1)$ by
action on $\cnj\PHone$}

Since $\cnj\PHone$ is a Hilbert space upon which the generators of
$su(1,1)$ $\{J_0,J_+,J_-\}$ act by left multiplication as unbounded
operators we can ask whether there exists a subspace of $\cnj\PHone$
for which they are bounded operators.

We propose the subspace $\Qexpr_\lambda$ given by
\begin{align}
\Qexpr_\lambda &= \spanrm\{Q^m_\lambda\ |\ m\in\Intg\}
\end{align}
where $\ad{J_0} Q^m_\lambda = \varepsilon m Q^m_\lambda$
and $\norm{Q^m_\lambda} =1 $.
Left multiplication by the generators of $su(1,1)$ on $\Qexpr_\lambda$
is given by
\begin{align}
J_0 Q^m_\lambda &= (\lambda+\varepsilon m) Q^m_\lambda &
J_\pm Q^m_\lambda &= (\lambda+\varepsilon m\pm \tfrac12 \pm i\Rp) 
Q^{m\pm1}_\lambda 
\label{rep_JQ}
\end{align}
where $\Rp^2=-R^2-\tfrac14\varepsilon^2\ge 0$.  These expressions are
similar to the standard continuous series of representation of
$su(1,1)$. 

The problem, which is still unsolved, is whether there exists
$\lambda\in\Cmpx$ for which $Q^m_\lambda$ has finite norm and hence
can by normalised.

By setting $Q^m_\lambda = \sum_{n=0}^\infty c_n
{P^m_n}/{\norm{P^m_n}}$, it is necessary to show that $|c_n|^2$
is a convergent series. 

Either by the manipulations of theorem \ref{thm_RM} or by manipulation
of the Hahn Polynomials one can show that
\begin{align}
J_0 P^m_n + P^m_n J_0 &= 
-\beta_{n+1}((n+1)^2-m^2)^\scrhalf P^m_{n+1}
-\wh\beta_{n}(n^2-m^2)^\scrhalf P^m_{n-1}
\label{SP_PJ0}
\end{align}
where
\begin{align*}
\beta_n &=
\frac{1}{2\alpha(2n)^\scrhalf(2n-1)^\scrhalf} \ ,
&\ 
\wh\beta_n &=
\frac{\alpha(4R^2+\varepsilon^2(1-n^2))}
{4(2n+1)(2n-1)^\scrhalf(2n)^\scrhalf}
&
\text{ and }
&&
\norm{P^m_n}^2\beta_n &= 
\norm{P^m_{n-1}}^2\wh\beta_n 
\end{align*}
After further manipulation we can show that the $c_n$
satisfy the recursive relation
\begin{align}
\gamma_{n+1} c_{n+1}
+i (\lambda + \tfrac12 \varepsilon m) c_n
+\gamma_{n} c_{n-1}
&= 0
\label{rep_reduction_form}
\end{align}
where
\begin{align*}
\gamma_n &= 
-i\beta_{n}(n^2-m^2)^\scrhalf 
\frac{\norm{P^m_n}}{\norm{P^m_{n-1}}} 
= \left(
\frac{(n^2-m^2)(4\Rp^2+\varepsilon^2 n^2)}
{16 (4n^2-1)}\right)^\scrhalf
\end{align*}
Substituting $c_n = n^a + O(n^{a-1})$ into (\ref{rep_reduction_form})
above we have $a=-\tfrac12+i(\lambda+m\varepsilon/2)$.  Thus the first
term in the expansion of $|c_n|^2$ is convergent if
$\Im(\lambda)>0$. This shows that the representation (\ref{rep_JQ})
cannot be a representation of the Lie group for which $\lambda$ must
be a real integer multiple of $\varepsilon$.

Further analysis is necessary to establish whether there is a
$\lambda\in\Cmpx$ for which  $|c_n|^2$ is
convergent series.


\section{Noncommutative Surfaces of Rotation}
\label{ch_suf}

As stated in the introduction, we would now like to consider what
other axially symmetric surfaces have noncommutative analogues. Here
we give a definition of an algebra $\Alg(\rho,\varepsilon)$ where
$\rho$ is an analytic function and $\varepsilon\in\Cmpx$, and
show that when $\varepsilon=0$ it is the commutative algebra of
functions on a surface of rotation. In subsection \ref{ch_suf_eg} we
give examples of the sphere, the hyperboloids, the paraboloid,
and the $q$-deformed sphere. We then show how to map between
noncommutative surfaces (subsection \ref{ch_suf_hom}), and whether
they have unitary representations (subsection
\ref{ch_suf_rep}). Finally we show how to interpret
$\Alg(\rho,\varepsilon)$ as the noncommutative complex plane, and give
noncommutative analogues of the stereographic projection of $S^2$ and
$\Hyptwop$ (subsections \ref{ch_suf_Cmplx} and \ref{ch_suf_Cmx_eg}).

\vskip 1em

Given an analytic function $\rho:\Cmpx\mapsto\Cmpx$ and a constant
$\varepsilon\in\Cmpx$ we define the algebra $\Alg(\rho,\varepsilon)$ to
be the set of polynomials generated by the elements 
\begin{align}
\{ X_0,X_+,X_-\}\union\{\rho(X_0+ r\varepsilon)\ |\ r\in\Intg\}
\label{suf_def_set}
\end{align}
quotiented by the ideal generated by
\begin{align}
[X_0,X_+] &= \varepsilon X_+ &
[X_0,X_-] &= -\varepsilon X_- &
X_+ X_- &= \rho(X_0) &
X_- X_+ &= \rho(X_0+\varepsilon) 
\label{suf_def_XX}
\end{align}
We say $\rho$ is {\it real} if $\rho|_\Real:\Real\mapsto\Real$. If
$\rho$ is real then there is a conjugation on $\Alg(\rho,\varepsilon)$
given by $X_0^\dagger=X_0,\, X_+^\dagger=X_-, \, X_-^\dagger=X_+$.
Also if $\rho$ is real, let $I_\rho\subset\Real$ be the set
$I_\rho=\{u\in\Real\ |\ \rho(u)>0\}$. This set is important for three
reasons: (1) When $\varepsilon=0$ it determines the topology of the
surface of rotation.  (2) It determines the nature of the unitary
representation of $\Alg(\rho,\varepsilon)$. (3) If $\rho|_{I_\rho}$ is
an invertible function then there exists an interpretation of
$\Alg(\rho,\varepsilon)$ in terms of noncommutative complex numbers.

If $I_\rho$ is connected then let $|I_\rho|$ be the size of $I_\rho$.
That is $|I_\rho|$ is the difference between the two endpoints if
$I_\rho$ is bounded and infinity otherwise.

\begin{theorem}
\label{th_suf_rot}
If $\rho$ is real, $I_\rho\ne\emptyset$ and $\varepsilon=0$ then
$\Alg(\rho,0)$ is the commutative algebra of polynomials in $(x,y,z)$
restricted to the surface 
\begin{align}
\Man_\rho &= \{(x,y,z)\in\Real\ |\ x^2+y^2=\rho(z) \}
\end{align}
where $X_0=z$ and $X_\pm=x\pm i y$. The limit of the commutator as
$\varepsilon\to0$ gives $\Man_\rho$ a Poison structure, given
by
\begin{align}
\{f,g\} &= 
\lim_{\varepsilon\to0} \left( \tfrac1\varepsilon [f,g] \right)
= 
i\left( 
\frac{\partial f}{\partial \phi} 
\frac{\partial g}{\partial z} 
-
\frac{\partial g}{\partial \phi} 
\frac{\partial f}{\partial z} 
\right)
\label{suf_Pois}
\end{align}
Furthermore, if $I_\rho$ is connected then one of the following
three is true:

$\bullet$ $I_\rho$ is bounded and $\Man_\rho$ is topologically
equivalent to the sphere

$\bullet$ $I_\rho$ is bounded only from one sided and $\Man_\rho$ is
topologically equivalent to the disc.

$\bullet$ $I_\rho=\Real$ and $\Man_\rho$ is
topologically equivalent to the cylinder.
\end{theorem}

\begin{proof}
From (\ref{suf_def_XX}), $\Alg(\rho,0)$ is a commutative algebra
and $X_+X_-= x^2+y^2 = \rho(z)$.

For (\ref{suf_Pois}) we note that both forms of the Poison bracket
are bi-differentials, that is they obey Leibniz rule with respect to
both variables. Therefore, it is only necessary to check the products of
the generators: $\{X_0,X_\pm\}$ and $\{X_+,X_-\}$.

The topology classes for $\Man_\rho$ are obvious.
\end{proof}

There exist more complicated situations if $I_\rho$ is not
connected. For instance, the surface may be locally topologically
equivalent to the intersections of two cones. These situations will
not be considered here.


\subsection{Examples}
\label{ch_suf_eg}

\subsubsection*{The Sphere and Hyperboloids}

We can see instantly that the noncommutative sphere $\PSph$ and
hyperboloids $\PHtwo$ and $\PHone$ are examples of
noncommutative surfaces with 
\begin{align}
\rho(u)=\alpha^2(R^2-u^2+\varepsilon u)
\label{suf_eg_S2}
\end{align}
where $\alpha^2,R^2,\varepsilon\in\Real$.

\subsubsection*{The paraboloid}

Let $\rho(u)=u$ then $\Man_\rho$ is a paraboloid.  From
(\ref{suf_def_XX}) we have $[X_+,X_-]=\varepsilon$ making $X_+$ and
$X_-$ the creation an annihilation operators for the Heisenberg-Weil
algebra.  Thus we can view the Heisenberg-Weil algebra as the
noncommutative paraboloid.

\subsubsection*{ The $q$ deformed Sphere, $su_q(2)$}

The algebra $su_q(2)$ is generated by $\{X_0,X_+,X_-\}$ which satisfy
\begin{align}
[X_0,X_\pm]&=\pm X_\pm &
[X_+,X_-]&=\frac{q^{2X_0} - q^{-2X_0}}{q - q^{-1}}
\end{align}
There are many ways of extending this to a set of
algebras, which are parameterised by $\varepsilon$ and which are
continuous when $\varepsilon=0$. One possibility is
\begin{align}
[X_0,X_\pm]&=\pm \varepsilon X_\pm &
[X_+,X_-]&=\frac{\sinh(\varepsilon\kappa)\sinh(2\kappa X_0)}
{\sinh(\kappa)^2}
\end{align}
where $e^\kappa=q$. To write this as a noncommutative surface of
rotation, let
\begin{align}
\rho(u) &= \frac{\cosh(3\kappa u-\varepsilon\kappa)}
{2\sinh(\kappa)^2}
+ \frac1{2\kappa^2} + R^2 + \frac{\varepsilon^2}4 
- \frac16 + C(\kappa,\varepsilon)
\end{align}
where $C(\kappa=0,\varepsilon)=0$. The constant (with respect to $u$)
in $\rho(u)$ is set by requiring that $\rho(u)\to R^2-u^2+\varepsilon
u$ as $\kappa\to 0$.

Setting $\varepsilon=0$, we have a deformed sphere for small $\kappa$,
whilst for large $\kappa$, $\Man_\rho$ tends to a cylinder (including
the discs at the top and bottom).  The cylinder has radius
$(R^2-\tfrac16+C(\kappa,0))^\scrhalf$, and length $2$.


\subsection{Homomorphism between noncommutative surfaces} 
\label{ch_suf_hom}

We give here a description for mapping between two noncommutative
surfaces of rotation.  These mappings are a generalisation of the
Holstein and Primakoff formalism \cite{Goldhirsh1}.

\begin{theorem}
\label{th_suf_hom}
Given algebras $\Alg(\rho_1,\varepsilon_1)$ generated by
$\{X_0,X_+,X_-\}$ and $\Alg(\rho_1,\varepsilon_2)$ generated by
$\{Y_0,Y_+,Y_-\}$ and given analytic functions
$\sigma_\pm:\Cmpx\mapsto\Cmpx\backslash \{0\}$ and $\lambda\in\Cmpx$
such that
\begin{align}
\rho_1\left(\frac{\varepsilon_1}{\varepsilon_2}u+\lambda\right) &= 
\rho_2(u) \sigma_+(u)\sigma_-(u)
\label{suf_hom_rrss}
\end{align}
there exists an homomorphism of algebras
\begin{align}
\Alg(\rho_1,\varepsilon)&\mapsto
\Alg(\rho_2,\sigma_+,\sigma_-,\varepsilon) 
&&&
\nonumber\\
X_0 &\mapsto \frac{\varepsilon_1}{\varepsilon_2}Y_0 + \lambda &
X_+ &\mapsto \sigma_+(Y_0)Y_+ &
X_- &\mapsto Y_-\sigma_-(Y_0)
\label{suf_hom_XY}
\end{align}
where
$\Alg(\rho_2,\sigma_+,\sigma_-,\varepsilon)$ is the enlarged algebra
generated by
\begin{align}
\Alg(\rho_2,\sigma_+,\sigma_-,\varepsilon)
&= \Alg(\rho_2,\varepsilon)\union
\{\sigma_+(Y_0+ r\varepsilon),\, \sigma_-(Y_0+ r\varepsilon)\ 
|\ r\in\Intg\}
\label{suf_hom_enlg_Ah2}
\end{align}

This mapping is injective but not necessarily surjective. If $\rho_1$
and $\rho_2$ are real then this mapping preserves conjugation if and
only if $\cnj{\sigma_+(u)}=\sigma_-(u)$ and $\lambda\in\Real$. If the
mapping preserves conjugation then $\Man_{\rho_1}$ is topologically equivalent to
$\Man_{\rho_2}$ and
\begin{align}
\frac{|I_{\rho_1}|}{\varepsilon_1} &=
\frac{|I_{\rho_2}|}{\varepsilon_2}
\end{align}
\end{theorem}  

\begin{proof}
This simply consists of substituting (\ref{suf_hom_XY}) into each
equation of (\ref{suf_def_XX}). Injectivity comes from the uniqueness
of polynomials. The topology come from looking at the zeros of $\rho$
which can only be shifted or rescaled. 
\end{proof}

In many cases we allow $\sigma_\pm$ to contain poles, zeros and branch
cuts (since they often contain a square root). This allows the mapping
from one topology to another, such as the stereographic projection of
the sphere and the bosonic representation of spin. The latter can be
viewed as a mapping between the paraboloid and the sphere.


\subsection{Representations of noncommutative surfaces; Crystals}
\label{ch_suf_rep}

For any algebra $\Alg(\rho,\varepsilon)$, there exists many non
unitary representations of this algebra: Given the functions 
$C,D:\Intg\mapsto\Cmpx$ such that $C(m)D(m)=\rho(\varepsilon m)$, then
a representation of $\Alg(\rho,\varepsilon)$ on the vector space
$\{\rvec(m)\}_{m\in\Intg}$ is given by
\begin{align}
X_0 \rvec(m) &= \varepsilon m \rvec(m) &
X_+ \rvec(m) &= C(m+1) \rvec(m+1) &
X_- \rvec(m) &= D(m) \rvec(m-1) 
\end{align}

The situation is more interesting if we wish our representation to be
unitary. The existence of such a representation implies that $\rho$ is
real and that $I_\rho$ is non empty. In the following we considered
only connected $I_\rho$.

\begin{theorem}
If $\rho|_\Real$ is real and $I_\rho$ is connected and non empty then
there exist a (unique up to phase) unitary representation of
$\Alg(\rho,\varepsilon)$ on the vector space $V$ with basis
$\{\rvec(m)\}_{m\in M}$, where $M\subset\Intg$.
This is given by
\begin{align}
X_0 \rvec(m) &= \varepsilon(m+\lambda) \rvec(m) &
X_+ \rvec(m) &= \cnj{D(m+1)} \rvec(m+1) &
X_- \rvec(m) &= D(m) \rvec(m-1) 
\label{suf_rep_rep}
\end{align}
where $\lambda\in\Real, 0\le\lambda<1$ and 
$|D(m)|^2=\rho(\varepsilon m+\varepsilon\lambda)$.
Only one of the following three must occur:

$\bullet$
If $|I_\rho|$ is finite then $M$ has is finite range
of $\Intg$ and $\varepsilon$ is constrained by
\begin{align}
|I_\rho|/\varepsilon
&= \dim V\in\Intg
\label{suf_rep_M}
\end{align}

$\bullet$
If $I_\rho$ is bounded from one side then so is $M$ and 
$\dim M=\infty$.

$\bullet$
If $I_\rho=\Real$ then $M\in\Intg$ and $\dim(V)=\infty$.
\end{theorem}

\begin{proof}
Clearly (\ref{suf_rep_rep}) is consistent with (\ref{suf_def_XX}).
Let $I_\rho$ be the range $-\infty\le u_{\min}\le u\le
u_{\max}\le\infty$ and $M$ be the range $-\infty\le m_{\min}\le m\le
m_{\max}\le\infty$.  Then from (\ref{suf_rep_rep}) we have
$\varepsilon m_{\min}=u_{\min}$ and $\varepsilon
(m_{\max}+1)=u_{\max}$.  So $m_{\min}$ (or $m_{\max}$) is finite if
$u_{\min}$ (or $u_{\max}$) is finite. If both are finite then
(\ref{suf_rep_M}) is obvious.
\end{proof}

When $\rho$ is given by (\ref{suf_eg_S2}) these representation
correspond to the standard unitary representation of the Lie algebra
$su(2)$ and $su(1,1)$. Further restrictions must be imposed to produce
the unitary representation of the Lie group $SU(1,1)$. 

\vskip 1em

There is a connection with one dimensional crystals of either finite
or infinite size. If the atoms are labelled by $\rvec(m)$ and the self
energy is proportional to $m$ (as could be the case for a simple
magnetic field) and the transition energy proportional to $D(m)$, then
we have a Hamiltonian of the form
\begin{align}
H &= X_0 + X_+ + X_-
\notag\\
&=
\sum_m \left(
\varepsilon(m+\lambda) \rvec(m)\lvec(m) +
\cnj{D(m+1)} \rvec(m+1)\lvec(m) +
D(m)\rvec(m-1)\lvec(m) \right)
\label{suf_crys}
\end{align}
This is an example of a combination of a Stark effect with a hopping
term.  On would like find the energy states for this Hamiltonian.
Clearly if $\rho=\alpha^2(R^2-u^2+\varepsilon u)$ then we can perform
a $su(2)$ or $su(1,1)$ rotation to produce the standard representation
of these groups. For general $\rho$ it may be possible to diagonalises
(\ref{suf_crys}) using first a Holstein-Primakoff transformation and
then the appropriate rotation.


\subsection{The noncommutative complex plane and stereographic
projections}
\label{ch_suf_Cmplx}

There is an alternative way of writing noncommutative surfaces such
that they look more like noncommutative domains in the complex plane.
Given the algebra $\Alg(\rho,\varepsilon)$ assume that $\rho$ is
invertible and that $\tau=\rho^{-1}$ then (\ref{suf_def_XX}) is
equivalent to the single equation
\begin{align}
\tau(z_-z_+) - \tau(z_+z_-) &= \varepsilon
\label{proj_tau_eps}
\end{align}
where $z_\pm=X_\pm$ and $X_0=\tau(z_+z_-)$.

If $\rho$ is real, then when $\varepsilon=0$ we reproduce the
commutative algebra of functions in $(z,\cnj z)$ on the domain
$\{|z|^2=\rho(u)\ \mbox{for some }u\in I_\rho\}\subset\Cmpx$. We can
see this by the substitution $z=z_-=e^{-i\phi}(\rho(z_0))^\scrhalf$,
$z_+=\cnj{z}$.

We can rewrite the projection given in theorem \ref{th_suf_hom}. This
projection is a noncommutative analogue of the stereographic
projection. Let $\Alg(\rho_1,\varepsilon)$ be a another surface of
rotation generated by $\{X_0,X_\pm\}$ and let
$\Alg(\tau^{-1},{\sigma_\pm},\varepsilon)$ be the extension of
$\Alg(\tau^{-1},\varepsilon)$ as before.
From (\ref{suf_hom_XY}) we have the mapping
$:\Alg(\rho_1,\varepsilon)\mapsto\Alg(\tau^{-1},
{\sigma_\pm},\varepsilon)$ given by
\begin{align}
X_+ &\mapsto \wt{\sigma_+}(z_+z_-) z_+ &
X_- &\mapsto z_- \wt{\sigma_-}(z_+z_-) &
X_0 &\mapsto \tau(z_+z_-)+\lambda
\end{align}
where $\rho_1(\tau(x))= x\wt{\sigma_+}(x)\wt{\sigma_-}(x)$ and
$\wt{\sigma_\pm}=\sigma_\pm\circ \tau$ for the functions $\sigma_\pm$
in theorem \ref{th_suf_hom}.


\subsection{Example: The Stereographic Projection of $\Sph$ and $\Hyptwop$}
\label{ch_suf_Cmx_eg}

We know that $\Alg(\rho,\varepsilon)$ with $\rho$ given by
(\ref{suf_eg_S2}), $R>0$ and $\varepsilon,\alpha^2\in\Real$ corresponds
to either $\PSph$ or $\PHtwo$ depending on the sign of $\alpha^2$.
The following map may be considered the noncommutative
analogue of a stereographic projection:
\begin{align}
J_0 &\mapsto \Rp\frac{4\Rp^2 - \alpha^2 x}{4\Rp^2 + \alpha^2 x} -
\frac{\varepsilon}2
&
J_+ &\mapsto
i
\frac{4\Rp^2\alpha^2}{4\Rp^2 + \alpha^2 x} z_+
&
J_- &\mapsto
-i z_-
\frac{4\Rp^2\alpha^2}{4\Rp^2 + \alpha^2 x} 
\end{align}
where $x=z_-z_+$ and $\Rp^2=R^2+\tfrac14\varepsilon^2$.
When $\varepsilon=0$ this map becomes the stereographic projection of
$S^2$ to $\Cmpx$ for $\alpha^2=1$ and the stereographic projection of
$\Hyptwop$ to the disc $\{|z|<2R\}\in\Cmpx$ for $\alpha^2=-1$.
For $\alpha^2=-1$ and $R^2\le0$ this is not a stereographic
projection.

Equation (\ref{proj_tau_eps}) in this case is equivalent to
\begin{align}
x - y &= 
\frac{-\varepsilon}{8\Rp^3\alpha^2}
(4\Rp^2 + \alpha^2 x)
(4\Rp^2 + \alpha^2 y)
\end{align}
where $x=z_-z_+$ and $y=z_+z_-$, or the M\"obius transformation
\begin{align}
y &= \frac{
(1+{\varepsilon}/{2\Rp})x + 2\varepsilon\Rp/\alpha^2}
{(-\varepsilon\alpha^2/8\Rp)x + (1-{\varepsilon}/{2\Rp})}  
\end{align}
We note that if $2\Rp=1$ and $\alpha^2=-1$ this is equivalent to the
algebra given in \cite{Klimek_Lesn1} and used latter in
\cite{Klimek_Lesn2} to give noncommutative version of surfaces with
higher genus.

The image of $P^m_n$ under this map may be written
\begin{align}
P^m_n &= 
\begin{cases}
(z_+)^m p^m_n(x) (\Rp^2+\alpha^2 x)^n & m\ge 0 \\
(z_-)^{-m} p^m_n(x) (\Rp^2+\alpha^2 x)^n & m< 0 
\end{cases}
\end{align}
where $p^m_n(x)$ is a polynomial of degree less than $n+1$, 
related to the Hahn polynomials.


\section{Discussion and Outlook}
\label{ch_Disc}

The most interesting case from section \ref{ch_Wick} is that of the
noncommutative one sheeted hyperboloid. This is a globally hyperbolic
spacetime and a two dimensional de Sitter space. As such it may be
related to inflation in the early universe. In order to do quantum
functional field theory (second quantisation) one must first construct
a Klein-Gordon inner product to distinguish positive and negative
frequency states. This may have the form
$\innerprod(f,g)=\pi_0(f^\dagger,\ad{J_0}g)$. One could then go on to
construct the Fock space. Clearly we would also like to finish the
calculation for the existence of unitary representation of $su(1,1)$
given by (\ref{rep_JQ}).

The new results about the noncommutative disc given as the image of
the two sheeted hyperboloid $\Hyptwop$ may given further insight of
higher genus surfaces using the analysis of Klimek and Lesniewski
\cite{Klimek_Lesn1,Klimek_Lesn2}

The product in $\PSph$ is equivalent to that discussed by Cahen
\cite{Cahen1}, who showed that was not a $\star$-product in the sense
of Flato et al. \cite{Flato1}. Since $\PHone$ and $\PHtwo$ are
algebraically equivalent to $\PSph$ these also cannot be
$\star$-product algebras. It would be useful to have an explicit
formula for this product in terms of an expansion in $\varepsilon$.

\vskip 1em

A principle objective of noncommutative geometry is the establishment
of a theory of quantum gravity. Starting from a noncommutative algebra
we would like to set up noncommutative analogues of concepts such as
vector fields, spinors, connections, curvature and ultimately
Einstein's equations and gravity.

Even deciding what is the analogue a vector field presents problems.
Vector fields have two properties 
which cannot both be required in 
noncommutative geometry: 

(1) they are derivatives of the algebra of
functions, and

(2) that form a module over the algebra of functions. 

\noindent For a matrix geometry all derivatives are inner and the
space of inner derivatives do not form a module of the algebra of
matrices.  This result is also true for the algebra $\PSph$.

Choosing vector fields to be derivatives \cite{Gratus4} then, for
matrix representations, there is a way of defining a three dimensional
space for 1-forms $\Omega^1(\Alg(\rho,\varepsilon))$. These are dual
to $\{\ad{X_0},\ad{X_+},\ad{X_-}\}$. However, \cite{Gratus4} also
shows that the dimension of the space of 2-forms
$\Omega^2(\Alg(\rho,\varepsilon))$ depends on the number of
``symmetries'' of the underlying space.  For the sphere and
hyperboloids, we can choose $\Omega^2(\PSph)$, $\Omega^2(\PHone)$ and
$\Omega^2(\PHtwo)$ to have up to four dimension, but in general
$\Omega^2(\Alg(\rho,\varepsilon))$ has at most two dimension.

Alternatively one could try and extend the approach of \cite{Gratus6}
and find ``fields'' which form a module over the algebra of functions but
which are derivatives only in the commutative limit.
An important step in this direction would be to establish a basis for the
algebra $\Alg(\rho,\varepsilon)$. One may start with functions of the
form $\{X_+^a X_0^b\,, X_-^aX_0^b\ |\ a,b\in\Intg^+\}$ but this set
would not include $\rho(X_0)$ unless it were a polynomial.
Also, one would like to establish which polynomials where
harmonic (like $P^m_n$ in section \ref{ch_Wick}). One would thus
generalise the  Laplace operator.  Its eigenstates would be
the harmonic (\ref{Pmn_Del_Pmn}), and should, for the finite
dimensional representations, also be orthogonal.

There is still no agreement on how to define connections, curvature,
etc. and there is much research in this area. However, having
noncommutative analogues of a large
collection of manifolds with different non constant
curvatures will enable one to examine many possible ideas.

Further problems will also be encountered when one wishes to construct
noncommutative analogues of spacetimes without a natural Poisson
structure. This includes the four dimensional spacetimes studied in
general relativity.  One might have to consider alternative approaches
such as adding additional dimensions or, more radically, considering
non-associative algebras.

\vskip 1em

As well as applications in the theory of crystals, noncommutative
surfaces of rotation may also have an interpretation in the theory of
strings, membranes, and higher $d$-branes.  The function
$\rho^\scrhalf$ might correspond to some kind of vibration on a closed
circular string which would not interact but may be created and
annihilated.

\subsection*{Acknowledgements}

The author would like to thank the Royal Society of London for a
European Junior Fellowship which enabled him to study in at Paris.
The author would also like to thank Robin Tucker, Peter
Pre\v{s}najder, Harald Grosse, and Karol Penson for useful discussions
which aided this work and, especially, Richard Kerner and the
Laboratoire de Gravitation et Cosmologie Relativistes, Paris~VI for
their hospitality.



\end{document}